\documentclass[10pt]{article}
\usepackage{epic,latexsym}
\usepackage{tikz}
\input epsf

\newtheorem{thm}{Theorem}
\newtheorem{prop}[thm]{Proposition}
\newtheorem{ob}[thm]{Observation}

\newtheorem{definition}{Definition}

\newcommand{\mod}{{\rm mod}}

\newcommand{\cC}{{\cal C}}
\newcommand{\cG}{{\cal G}}

\newcommand{\cH}{{\cal H}}
\newcommand{\cF}{{\cal F}}
\newcommand{\cL}{{\cal L}}

\newcommand{\gt}{\gamma_t}
\newcommand{\gta}{\gamma_t^a}

\newenvironment{unnumbered}[1]{\trivlist
\item [\hskip \labelsep {\bf
#1}]\ignorespaces\it}{\endtrivlist}

\newcommand{\1}{\vspace{0.1cm}}
\newcommand{\2}{\vspace{0.2cm}}

\newcommand{\smallqed}{{\tiny ($\Box$)}}
\newcommand{\qed}{$\Box$}

\newcommand{\cub}{{\rm cubic}}
\newcommand{\done}{\delta = 1}
\newcommand{\dtwo}{\delta = 2}

\def\vertex(#1){\put(#1){\circle*{4}}}
\def\vertexo(#1){\put(#1){\circle{2}}}
\def\vert(#1){\put(#1){\circle*{1.5}}}
\def\verto(#1){\put(#1){\circle{1.5}}}
\def\lab(#1)#2{\put(#1){\makebox(0,0)[c]{#2}}}
\begin{document}

\bibliographystyle{plain}

\title{Equality in a Linear Vizing-Like Relation that Relates the \\ Size and
Total Domination Number of a Graph}

\author{Michael A. Henning\thanks{Research supported in part
by the University of Johannesburg and the South African National
Research Foundation} \, and Ernst  J. Joubert \\
Department of Mathematics \\
University of Johannesburg \\
Auckland Park 2006, South Africa }
\date{}
\maketitle

\begin{abstract}
Let $G$ be a graph each component of which has order at least~$3$,
and let $G$ have order $n$, size $m$, total domination number $\gt$
and maximum degree $\Delta(G)$. Let $\Delta = 3$ if $\Delta(G) = 2$
and $\Delta = \Delta (G)$ if $\Delta(G) \ge 3$. It is known [J. Graph
Theory 49 (2005), 285--290; J. Graph Theory 54 (2007), 350--353] that
$m \le \Delta (n- \gt)$. In this paper we characterize the extremal
graphs $G$ satisfying $m = \Delta (n-\gt)$.
\end{abstract}

{\small \textbf{Keywords:} Maximum degree; order; size; total domination.} \\
\indent {\small \textbf{AMS subject classification: 05C69}}

\section{Introduction}

In this paper we continue the study of total domination in graphs.
Let $G = (V,E)$ be a graph with vertex set $V$, edge set $E$ and no
isolated vertex. A \emph{total dominating set}, abbreviated TD-set,
of $G$ is a set $S$ of vertices of $G$ such that every vertex is
adjacent to a vertex in $S$. The \emph{total domination number} of
$G$, denoted by $\gt(G)$, is the minimum cardinality of a TD-set. A
TD-set of $G$ of cardinality $\gt(G)$ is called a $\gt(G)$-set.
Total domination in graphs is now well studied in graph theory. The
literature on this subject has been surveyed and detailed in the two
books by Haynes, Hedetniemi, and Slater \cite{hhs1,hhs2}. A recent
survey of total domination in graphs can be found in~\cite{He09}.

A classical result of Vizing~\cite{Vi} relates the size and the
ordinary domination number, $\gamma$, of a graph of given order.
Rautenbach~\cite{Ra90} shows that the square dependence on $n$ and
$\gamma$ in the result of Vizing turns into a linear dependence on
$n$, $\gamma$, and the maximum degree $\Delta$.

Dankelmann et al.~\cite{ddgghs} proved a Vizing-like relation between
the size and the total domination number of a graph of given order.
Sanchis~\cite{Sa04} showed that if we restrict our attention to
connected graphs with total domination number at least~$5$, then the
bound in~\cite{ddgghs} can be improved slightly. The square
dependence on $n$ and $\gt$ presented in~\cite{ddgghs,Sa04} is
improved in~\cite{He05,ShKaHe07,Ye07} into a linear dependence on
$n$, $\gt$ and $\Delta$ by demanding a more even distribution of the
edges by restricting the maximum degree $\Delta$. In particular, the
following linear Vizing-like relation relating the size of a graph
and its order, total domination number, and maximum degree is
established in~\cite{He05,ShKaHe07}.

\begin{unnumbered}{Theorem~A (\cite{He05,ShKaHe07})}
Let $G$ be a graph each component of which has order at least~$3$,
and let $G$ have order~$n$, size~$m$, total domination number $\gt$,
and maximum degree~$\Delta(G)$. Let $\Delta = 3$ if $\Delta(G) = 2$
and $\Delta = \Delta (G)$ if $\Delta(G) \ge 3$. Then, $m \le \Delta
(n- \gt)$.
\end{unnumbered}

Our aim in this paper is to characterize the extremal graphs
achieving equality in the upper bound in Theorem~A; that is, to
characterize the graphs $G$ satisfying the statement of Theorem~A
such that $m =\Delta (n - \gt)$.

\subsection{Notation}

For notation and graph theory terminology we in general
follow~\cite{hhs1}. Specifically, let $G = (V,E)$ be a graph with
vertex set $V$ of order~$n(G) = |V|$ and edge set $E$ of size~$m(G) =
|E|$, and let $v$ be a vertex in $V$.
The \emph{open neighborhood} of $v$ is $N_G(v) = \{u \in V \, | \, uv
\in E(G)\}$ and the \emph{closed neighborhood of $v$} is $N_G[v] =
\{v\} \cup N(v)$. The degree of $v$ is $d_G(v) = |N_G(v)|$. The
minimum and maximum degree among the vertices of $G$ is denoted by
$\delta(G)$ and $\Delta(G)$, respectively.
A vertex adjacent to a vertex of degree~$1$ is called a
\emph{support} vertex.
For a set $S \subseteq V$, its \emph{open neighborhood} is the set
$N_G(S) = \cup_{v \in S} N_G(v)$, and its \emph{closed neighborhood}
is the set $N_G[S] = N_G(S) \cup S$.
If the graph $G$ is clear from the context, we simply write $N(v)$
and $d(v)$ rather than $N_G(v)$ and $d_G(v)$, respectively. Further
we write $N[v]$, $N[S]$ and $N(S)$ rather that $N_G[v]$, $N_G[S]$ and
$N_G(S)$, respectively.
For sets $A, B \subseteq V$, we say that $A$ \emph{dominates} $B$ if
$B \subseteq N[A]$, while $A$ \emph{totally dominates} $B$ if $B
\subseteq N(A)$.

For a set $S \subseteq V$, the subgraph induced by $S$ is denoted by
$G[S]$. Further if $S \ne V$, then we denote the graph obtained from
$G$ by deleting all vertices in $S$ by $G - S$. A component of $G$
that is isomorphic to a graph $F$ is called an $F$-component of $G$.

A \emph{cycle} on $n$ vertices is denoted by $C_n$, while a
\emph{path} on $n$ vertices is denoted by $P_n$. We denote by $K_n$
the \emph{complete graph} on $n$ vertices.
A $2$-\emph{path} in $G$ is a path on at least three vertices with
both ends of the path having degree at least~$3$ in $G$ and with
every internal vertex of the path having degree~$2$ in~$G$.
A \emph{special} $2$-\emph{path} in $G$ is a $2$-path
$v_1v_2v_3v_4v_5$ such that $v_1$ and $v_5$ have two common
neighbors, $x$ and $y$ say, in $G$, and the vertices $v_1$, $v_5$,
$x$ and $y$ all have degree~$3$ in $G$. In particular, we note that
$N(v_1) = \{v_2,x,y\}$ and $N(v_5) = \{v_4,x,y\}$.

\section{Special Graphs and Families of Graphs}

\subsection{The Family~$\cG_\cub$}
\label{S:GH}

Let $GP_{16}$ denote the generalized Petersen graph of order~$16$
shown in Figure~\ref{GenPeter}.

\setlength{\unitlength}{0.7mm}
\begin{figure}[htb]
\begin{center}

\begin{picture}(0,43)


\multiput(-20,10)(0,20){2}{\circle*2}
\multiput(-10,0)(0,40){2}{\circle*2}
\multiput(-10,15)(0,10){2}{\circle*2}
\multiput(-5,10)(0,20){2}{\circle*2}

\multiput(20,10)(0,20){2}{\circle*2}
\multiput(10,0)(0,40){2}{\circle*2}
\multiput(10,15)(0,10){2}{\circle*2}
\multiput(5,10)(0,20){2}{\circle*2}

\put(-20,10){\line(1,-1){10}} \put(-20,10){\line(2,1){10}}
\put(-20,10){\line(0,1){20}} \put(-20,30){\line(1,1){10}}
\put(-20,30){\line(2,-1){10}} \put(-10,0){\line(1,0){20}}
\put(-10,40){\line(1,0){20}} \put(-10,0){\line(1,2){5}}
\put(-10,40){\line(1,-2){5}}

\put(-10,15){\line(1,0){20}} \put(-10,25){\line(1,0){20}}
\put(-5,10){\line(0,1){20}} \put(5,10){\line(0,1){20}}

\put(-10,15){\line(1,1){15}} \put(-10,25){\line(1,-1){15}}
\put(-5,10){\line(1,1){15}} \put(-5,30){\line(1,-1){15}}

\put(20,10){\line(-1,-1){10}} \put(20,10){\line(-2,1){10}}
\put(20,10){\line(0,1){20}} \put(20,30){\line(-1,1){10}}
\put(20,30){\line(-2,-1){10}} \put(10,0){\line(-1,2){5}}
\put(10,40){\line(-1,-2){5}}

\end{picture}


\end{center}
\vskip -0.3 cm \caption{The generalized Petersen graph $GP_{16}$ of
order~$16$.} \label{GenPeter}
\end{figure}
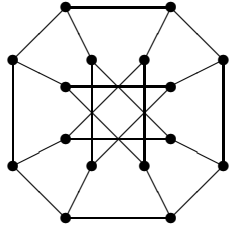

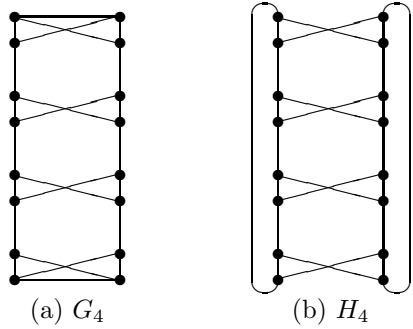
\begin{figure}[htb]
\begin{center}
\setlength{\unitlength}{0.35mm}
\begin{picture}(0,105)


\multiput(-70,0)(0,30){4}{\circle*4}
\multiput(-70,10)(0,30){4}{\circle*4}

\multiput(-30,0)(0,30){4}{\circle*4}
\multiput(-30,10)(0,30){4}{\circle*4}


\put(-70,1){\line(0,1){8}} \put(-70,11){\line(0,1){18}}
\put(-70,31){\line(0,1){24}} \put(-70,41){\line(0,1){6}}
\put(-70,53){\line(0,1){6}} \put(-70,61){\line(0,1){8}}
\put(-70,71){\line(0,1){18}} \put(-70,91){\line(0,1){8}}

\put(-30,1){\line(0,1){8}} \put(-30,11){\line(0,1){18}}
\put(-30,31){\line(0,1){24}}
\put(-30,41){\line(0,1){6}} \put(-30,53){\line(0,1){6}}
\put(-30,61){\line(0,1){8}} \put(-30,71){\line(0,1){18}}
\put(-30,91){\line(0,1){8}}

\put(-69,0.5){\line(4,1){38}} \put(-69,9.5){\line(4,-1){38}}
\put(-69,30.5){\line(4,1){38}} \put(-69,39.5){\line(4,-1){38}}
\put(-69,60.5){\line(4,1){38}} \put(-69,69.5){\line(4,-1){38}}
\put(-69,90.5){\line(4,1){38}} \put(-69,99.5){\line(4,-1){38}}

\put(-69,0){\line(1,0){38}} \put(-69,100){\line(1,0){38}}

\lab(-50,-12){(a) $G_4$}


\multiput(30,0)(0,30){4}{\circle*4}
\multiput(30,10)(0,30){4}{\circle*4}

\multiput(70,0)(0,30){4}{\circle*4}
\multiput(70,10)(0,30){4}{\circle*4}


\put(30,1){\line(0,1){8}} \put(30,11){\line(0,1){18}}
\put(30,31){\line(0,1){24}} \put(30,41){\line(0,1){6}}
\put(30,53){\line(0,1){6}} \put(30,61){\line(0,1){8}}
\put(30,71){\line(0,1){18}} \put(30,91){\line(0,1){8}}

\put(70,1){\line(0,1){8}} \put(70,11){\line(0,1){18}}
\put(70,31){\line(0,1){24}}
\put(70,41){\line(0,1){6}} \put(70,53){\line(0,1){6}}
\put(70,61){\line(0,1){8}} \put(70,71){\line(0,1){18}}
\put(70,91){\line(0,1){8}}

\put(31,0.5){\line(4,1){38}} \put(31,9.5){\line(4,-1){38}}
\put(31,30.5){\line(4,1){38}} \put(31,39.5){\line(4,-1){38}}
\put(31,60.5){\line(4,1){38}} \put(31,69.5){\line(4,-1){38}}
\put(31,90.5){\line(4,1){38}} \put(31,99.5){\line(4,-1){38}}

\put(25,100){\oval(10,10)[tr]} \put(25,100){\oval(10,10)[tl]}
\put(75,100){\oval(10,10)[tr]} \put(75,100){\oval(10,10)[tl]}

\put(25,0){\oval(10,10)[br]} \put(25,0){\oval(10,10)[bl]}
\put(75,0){\oval(10,10)[br]} \put(75,0){\oval(10,10)[bl]}

\put(20,0){\line(0,1){100}} \put(80,0){\line(0,1){100}}

\lab(50,-12){(b) $H_4$}


\end{picture}


\end{center}
\caption{Cubic graphs $G_4 \in \cG$ and $H_4 \in \cH$.} \label{GH}
\end{figure}

The following two infinite families $\cG$ and $\cH$ of connected
cubic graphs (described below) with total domination number one-half
their orders are constructed in~\cite{fhmp}. For $k \ge 1$, let $G_k$
be the graph constructed as follows. Consider two copies of the path
$P_{2k}$ with respective vertex sequences $a_{1}b_{1}a_{2}b_{2}
\ldots a_{k}b_{k}$ and $c_{1}d_{1}c_{2}d_{2} \ldots c_{k}d_{k}$. Let
$A = \{a_1,a_2,\ldots,a_k\}$, $B = \{b_1,b_2,\ldots,b_k\}$, $C =
\{c_1,c_2,\ldots,c_k\}$, and $D = \{d_1,d_2,\ldots,d_k\}$. For each
$i \in \{1,2,\ldots,k\}$, join $a_{i}$ to $d_{i}$ and $b_{i}$ to
$c_{i}$. To complete the construction of the graph $G_k \in \cG$ join
$a_{1}$ to $c_{1}$ and $b_{k}$ to $d_{k}$. Let $\cG = \{ G_k \mid k
\ge 1\}$. For $k \ge 2$, let $H_k$ be obtained from $G_k$ by deleting
the two edges $a_{1}c_{1}$ and $b_{k}d_{k}$ and adding the two edges
$a_{1}b_{k}$ and $c_{1}d_{k}$. Let $\cH = \{ H_k \mid k \ge 2\}$. We
note that $G_k$ and $H_k$ are cubic graphs of order~$4k$. Further, we
note that $G_1 = K_4$. The graphs $G_4 \in \cG$ and $H_4 \in \cH$,
for example, are illustrated in Figure~\ref{GH}.

Let $\cG_\cub = \cG \cup \cH \cup \{GP_{16}\}$. We note that each
graph in the family $\cG_\cub$ is a cubic graph.

\subsection{The Family~$\cG_{\dtwo}$}
\label{S:Gtwo}

By \emph{contracting} two vertices $x$ and $y$ in $G$ we mean
replacing the vertices $x$ and $y$ by a new vertex $v_{xy}$ and
joining $v_{xy}$ to all vertices in $V(G) \setminus \{x,y\}$ that
were adjacent to $x$ or $y$ in $G$. Let $\cG_3$ be a set of graphs
only containing one element, namely the $3$-cycle $C_3$. Similarly,
let $\cG_6 = \{C_6\}$. For notational convenience, let $\cG_4 =
\emptyset$ and let $\cG_5 = \emptyset$. For every $i>6$, define
$\cG_i$ as follows.\footnote{We remark that for $i \ge 3$, our family
$\cG_i$ is a subfamily of the family called $\cC_i$ constructed in
Section~3 in~\cite{HeYe09}.}

\begin{definition} \label{defn:cCi}
For every $i>6$, a graph $R_i$ belongs to $\cG_i$ if and only if
$\delta(R_i) \ge 2$ and $R_i$ contains a special $2$-path
$v_1v_2v_3v_4v_5$ and the graph obtained by contracting $v_1$ and
$v_5$ in $R_i$ and deleting $\{v_2,v_3,v_4\}$ belongs to $\cG_{i-4}$.
\end{definition}

Let $i \ge 3$. We note that $\cG_i =  \emptyset$ for $i \equiv 0,1 \,
(\mod \, 4)$. For $i \equiv 2,3 \, (\mod \, 4)$, suppose $R_i$
belongs to $\cG_i$. Let $v_1v_2v_3v_4v_5$ be a special $2$-path in
$R_i$ and let $R_{i-4}$ be the graph in $\cG_{i-4}$ obtained by
contracting $v_1$ and $v_5$ in $R_i$ into a new vertex $w$ and
deleting $\{v_2,v_3,v_4\}$. Further let $x$ and $y$ be the two common
neighbors of $v_1$ and $v_5$ in $R_i$. We note then that $x$ and $y$
both have degree~$2$ in $R_{i-4}$ and have $w$ as a common neighbor.

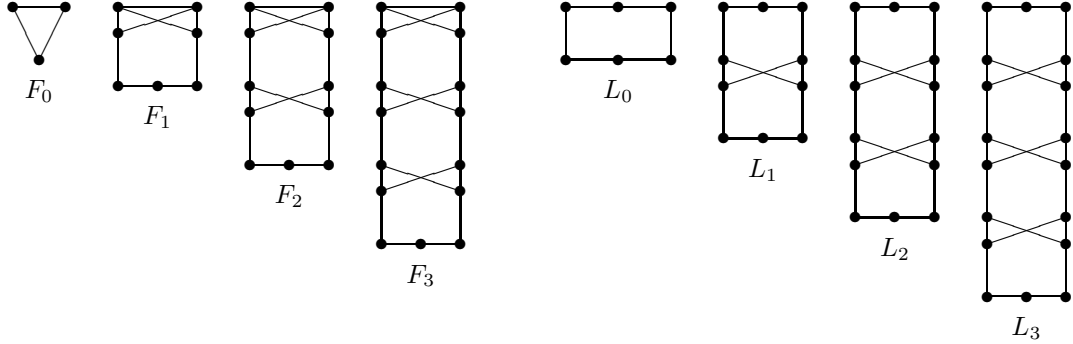
\begin{figure}[htb]
\begin{center}
\setlength{\unitlength}{0.35mm}
\begin{picture}(10,107)


\multiput(-190,100)(20,0){2}{\circle*4} \vertex(-180,80)
\put(-190,100){\line(1,0){20}} \put(-180,80){\line(-1,2){10}}
\put(-180,80){\line(1,2){10}}

\lab(-180,68){$F_0$}


\multiput(-150,90)(0,30){1}{\circle*4}
\multiput(-150,70)(0,30){2}{\circle*4} \vertex(-135,70)

\multiput(-120,90)(0,30){1}{\circle*4}
\multiput(-120,70)(0,30){2}{\circle*4}

\put(-150,70){\line(0,1){30}} \put(-120,70){\line(0,1){30}}

\put(-150,90){\line(3,1){30}} \put(-150,100){\line(3,-1){30}}

\put(-150,70){\line(1,0){30}} \put(-150,100){\line(1,0){30}}

\put(-150,70){\line(0,1){30}} \put(-120,70){\line(0,1){30}}

\lab(-135,58){$F_1$}



\multiput(-100,60)(0,30){2}{\circle*4}
\multiput(-100,40)(0,30){3}{\circle*4} \vertex(-85,40)

\multiput(-70,60)(0,30){2}{\circle*4}
\multiput(-70,40)(0,30){3}{\circle*4}

\put(-100,40){\line(0,1){60}} \put(-70,40){\line(0,1){60}}

\put(-70,41){\line(0,1){6}} \put(-70,53){\line(0,1){6}}
\put(-70,61){\line(0,1){8}} \put(-70,71){\line(0,1){18}}
\put(-70,91){\line(0,1){8}}

\put(-100,60){\line(3,1){30}} \put(-100,70){\line(3,-1){30}}
\put(-100,90){\line(3,1){30}} \put(-100,100){\line(3,-1){30}}

\put(-100,40){\line(1,0){30}} \put(-100,100){\line(1,0){30}}

\lab(-85,28){$F_2$}



\multiput(-50,30)(0,30){3}{\circle*4}
\multiput(-50,10)(0,30){4}{\circle*4} \vertex(-35,10)

\multiput(-20,30)(0,30){3}{\circle*4}
\multiput(-20,10)(0,30){4}{\circle*4}

\put(-50,11){\line(0,1){18}} \put(-50,31){\line(0,1){24}}
\put(-50,41){\line(0,1){6}} \put(-50,53){\line(0,1){6}}
\put(-50,61){\line(0,1){8}} \put(-50,71){\line(0,1){18}}
\put(-50,91){\line(0,1){8}}

\put(-20,11){\line(0,1){18}} \put(-20,31){\line(0,1){24}}
\put(-20,41){\line(0,1){6}} \put(-20,53){\line(0,1){6}}
\put(-20,61){\line(0,1){8}} \put(-20,71){\line(0,1){18}}
\put(-20,91){\line(0,1){8}}

\put(-50,30){\line(3,1){30}} \put(-50,40){\line(3,-1){30}}
\put(-50,60){\line(3,1){30}} \put(-50,70){\line(3,-1){30}}
\put(-50,90){\line(3,1){30}} \put(-50,100){\line(3,-1){30}}

\put(-50,10){\line(1,0){30}} \put(-50,100){\line(1,0){30}}

\lab(-35,-2){$F_3$}



\multiput(20,100)(20,0){3}{\circle*4}
\multiput(20,80)(20,0){3}{\circle*4}

\put(20,100){\line(1,0){40}} \put(20,80){\line(1,0){40}}
\put(20,80){\line(0,1){20}} \put(60,80){\line(0,1){20}}

\lab(40,68){$L_0$}



\multiput(80,100)(15,0){3}{\circle*4}
\multiput(80,80)(30,0){2}{\circle*4}
\multiput(80,70)(30,0){2}{\circle*4}
\multiput(80,50)(15,0){3}{\circle*4}

\put(80,100){\line(1,0){30}} \put(80,50){\line(1,0){30}}
\put(80,50){\line(0,1){50}} \put(110,50){\line(0,1){50}}

\put(80,70){\line(3,1){30}} \put(80,80){\line(3,-1){30}}

\lab(95,38){$L_1$}


\multiput(130,100)(15,0){3}{\circle*4}
\multiput(130,80)(30,0){2}{\circle*4}
\multiput(130,70)(30,0){2}{\circle*4}

\multiput(130,50)(30,0){2}{\circle*4}
\multiput(130,40)(30,0){2}{\circle*4}

\multiput(130,20)(15,0){3}{\circle*4}

\put(130,100){\line(1,0){30}} \put(130,20){\line(1,0){30}}
\put(130,20){\line(0,1){80}} \put(160,20){\line(0,1){80}}

\put(130,70){\line(3,1){30}} \put(130,80){\line(3,-1){30}}
\put(130,40){\line(3,1){30}} \put(130,50){\line(3,-1){30}}

\lab(145,8){$L_2$}


\multiput(180,100)(15,0){3}{\circle*4}
\multiput(180,80)(30,0){2}{\circle*4}
\multiput(180,70)(30,0){2}{\circle*4}

\multiput(180,50)(30,0){2}{\circle*4}
\multiput(180,40)(30,0){2}{\circle*4}

\multiput(180,20)(30,0){2}{\circle*4}
\multiput(180,10)(30,0){2}{\circle*4}

\multiput(180,-10)(15,0){3}{\circle*4}

\put(180,100){\line(1,0){30}} \put(180,-10){\line(1,0){30}}
\put(180,-10){\line(0,1){110}} \put(210,-10){\line(0,1){110}}

\put(180,70){\line(3,1){30}} \put(180,80){\line(3,-1){30}}
\put(180,40){\line(3,1){30}} \put(180,50){\line(3,-1){30}}
\put(180,10){\line(3,1){30}} \put(180,20){\line(3,-1){30}}

\lab(195,-22){$L_3$}

\end{picture}


\end{center}
\vskip 0.5 cm  \caption{The graphs $F_0,F_1,F_2,F_3$ and
$L_0,L_1,L_2,L_3$.} \label{F0H0}
\end{figure}

For each $i \ge 0$, the family $\cG_{4i+3}$ consists of precisely one
graph $F_i$, namely the graph $F_i = C_3$ when $i = 0$ and, for $i
\ge 1$, the graph $F_i$ which is obtained from the graph $G_i$
defined in Section~\ref{S:GH} by subdividing the edge $a_1c_1$ three
times. We also note that for each $i \ge 0$, the family $\cG_{4i+2}$
consists of precisely one graph $L_i$, namely the graph $L_i = C_6$
when $i = 0$ and, for $i \ge 1$, the graph $L_i$ which is obtained
from the graph $G_i$ defined in Section~\ref{S:GH} by subdividing the
edge $a_1c_1$ three times and subdividing the edge $b_kd_k$ three
times.
The graphs $F_0,F_1,F_2,F_3$ and $L_0,L_1,L_2,L_3$, for example, are
shown in Figure~\ref{F0H0}.

Let $\cF = \{ F_k \mid k \ge 0\}$ and let $\cL = \{ L_k \mid k \ge
0\}$. Let $\cG_{\dtwo} = \cF \cup \cL$. We note that each graph in
the family $\cG_{\dtwo}$ has minimum degree~$\delta = 2$.

\subsection{The Family~$\cG_{\done}$}

For a graph $H$, we denote by $H \circ P_2$ the graph of
order~$3|V(H)|$ obtained from $H$ by attaching a path of length~$2$
to each vertex of $H$ so that the resulting paths are
vertex-disjoint. The graph $H \circ P_2$ is also called the
$2$-\emph{corona} of $H$. The graph $C_4 \circ P_2$ is shown in
Figure~\ref{f:C4P2}.

\setlength{\unitlength}{0.6mm}
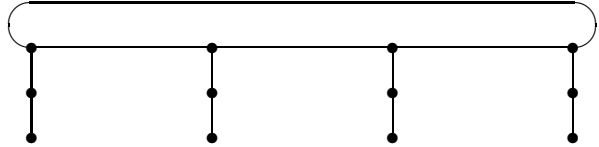
\begin{figure}[htb]
\begin{center}
\begin{picture}(0,35)
\multiput(-60,20)(40,0){4}{\circle*{2.5}}
\multiput(-60,10)(40,0){4}{\circle*{2.5}}
\multiput(-60,0)(40,0){4}{\circle*{2.5}}

\put(-60,20){\line(0,-1){20}} \put(-20,20){\line(0,-1){20}}
\put(60,20){\line(0,-1){20}} \put(20,20){\line(0,-1){20}}

\put(-60,25){\oval(10,10)[bl]} \put(-60,25){\oval(10,10)[tl]}
\put(60,25){\oval(10,10)[br]} \put(60,25){\oval(10,10)[tr]}
\put(-60,30){\line(1,0){120}} \put(-60,20){\line(1,0){120}}

\end{picture}
\end{center}
\vskip -.3 cm \caption{The graph $C_4 \circ P_2$.} \label{f:C4P2}
\end{figure}

Let $\cG_{\done}$ denote the family of all graphs $H \circ P_2$,
where $H$ is a cycle $C_k$ on $k \ge 3$ vertices. We note that each
graph in the family $\cG_{\done}$ has minimum degree~$\delta = 1$.

\section{Main Result}

We shall prove the following result, a proof of which is presented in
Section~\ref{S:mainp}.

\begin{thm}
Let $G$ be a connected graph of order~$n$, size~$m$, total domination
number $\gt$, and maximum degree~$\Delta(G)$ with each component of
$G$ of order at least~$3$. Let $\Delta = 3$ if $\Delta(G) = 2$ and
$\Delta = \Delta (G)$ if $\Delta(G) \ge 3$. Then, $m \le \Delta (n-
\gt)$, with equality if and only if $G \in \cG_{\done} \cup
\cG_{\dtwo} \cup \cG_{\cub}$. \label{thm1}
\end{thm}

\section{Known Results and Preliminary Results}

The total domination number of a path $P_n$ and a cycle $C_n$ on $n
\ge 3$ vertices is easy to compute (or see, \cite{He00}): For $n \ge
3$, $\gt(P_n) = \gt(C_n) = \lfloor n/2 \rfloor + \lceil n/4 \rceil -
\lfloor n/4 \rfloor$. We shall need the following upper bounds on the
total domination number of a graph.

\begin{unnumbered}{Theorem~B}
Let $G$ be a graph of order $n$ and let $F$ be the subgraph of $G$
induced by its vertices of degree~$2$, if such vertices exist. Then
the following holds. \\
\indent {\rm (a) (\cite{cdh})} If every component of
$G$ has order at least~$3$, then $\gt(G) \le 2n/3$. \\
\indent {\rm (b) (\cite{Alfewy,ChMc,Tuza})} If $\delta(G) \ge 3$,
then $\gt(G) \le n/2$. \\
\indent {\rm (c) (\cite{HeYe08})} If $G$ is a connected graph with
$\delta(G) \ge 3$, then $\gt(G) = n/2$ if and only if $G \in
\cG_\cub$. \\
\indent {\rm (d) (\cite{LaWe})} If $\delta(G) \ge 2$ and every
component of $F$ has order at most~$2$, then $\gt(G) \le n/2$.
\end{unnumbered}

For a graph $G = (V,E)$ and a vertex $v \in V$, an \emph{almost total
dominating set} of $G$, abbreviated ATD-set, with respect to $v$ is a
set $S$ of vertices of $G$ such that $v \in S$ and every vertex
different from $v$ is adjacent to a vertex in $S$ while $v$ is
isolated in $G[S]$. The \emph{almost total domination number} of $G$
with respect to $v$, denoted by $\gta(G;v)$, is the minimum
cardinality of an ATD-set with respect to $v$. An ATD-set of $G$ with
respect to $v$ of cardinality $\gta(G;v)$ is called a
$\gta(G;v)$-set. Every $\gta(G;v)$-set can be extended to a TD-set of
$G$ by adding to it a neighbor of $v$, and so $\gt(G) \le \gta(G;v) +
1$. We shall also need the following properties of graphs in the
families~$\cG_{\done} \cup \cG_{\dtwo}$.

\begin{ob}
Let $G \in \cG_{\done}$ have order $n$ and let $v \in V(G)$. Then the
following holds. \\
\indent {\rm (a)} If $G \in \cG_{\done}$, then $\gt(G) = 2n/3$. \\
\indent {\rm (b)} If $d_G(v) = 1$, then $\gt(G - v) = \gt(G) - 1$. \\
\indent {\rm (c)} If $d_G(v) = 2$, then $\gta(G;v) = \gt(G) - 1$.
 \label{ob:Gone}
\end{ob}

We note that if $G \in \cF$ and $G \ne C_3$, then by construction
there is only one $2$-path in $G$ and this $2$-path is a special
$2$-path.

\begin{prop}
Let $G \in \cF$ have order~$n$. If $G \ne C_3$, let $uv_1v_2v_3v$ be
a special $2$-path in $G$ and let $w \in \{v_1,v_2,v_3\}$. Then
the following holds. \\
\indent {\rm (a)} $\gt(G) = (n+1)/2$. \\
\indent {\rm (b)} $\gt(G - w) = (n-1)/2$. \\
\indent {\rm (c)} $\gta(G;w) = (n-1)/2$. \label{p:cF}
\end{prop}
\textbf{Proof.} Suppose that $G \in \cF$ has order~$n$. Then, $G =
F_k$ for some integer $k \ge 0$ and $n = 4k+3$. If $k = 0$, then $G =
C_3$. If $k \ge 1$, then following the notation introduced earlier,
the graph $F_k$ is obtained from the graph $G_k$ defined in
Section~\ref{S:GH} by subdividing the edge $a_1c_1$ three times. Let
$a_1v_1v_2v_3c_1$ denote the resulting path in $F_k$. We note that
this is the only $2$-path in $G$ and this $2$-path is in fact a
special $2$-path in $G$.

(a) If $k = 0$, then $G = C_3$ and $\gt(G) = 2 = (n+1)/2$, as
desired. Hence we may assume that $k \ge 1$. We show first that
$\gt(F_k) \le \gt(G_k) + 2$. Let $S$ be a $\gt(G_k)$-set. If
$\{a_1,c_1\} \subseteq S$, then we can simply replace $c_1$ in $S$
with the vertex $b_1$ (or $d_1$). Hence we may choose $S$ to contain
at most one of $a_1$ and $c_1$. But then $S \cup \{v_1,v_2\}$ is a
TD-set of $F_k$, and so $\gt(F_k) \le |S| + 2 = \gt(G_k) + 2$.
We show next that $\gt(F_k) \ge \gt(G_k) + 2$. Among all
$\gt(F_k)$-sets, let $D$ be chosen to contain as few vertices from
the set $\{v_1,v_2,v_3\}$ as possible. In order to totally dominate
$v_2$, we may assume, renaming vertices if necessary, that $v_1 \in
D$. If $v_3 \in D$, then replacing $v_3$ in $D$ with $b_1$ (or $d_1$)
produces a new $\gt(F_k)$-set that contradicts our choice of $D$.
Hence, $v_3 \notin D$. In order to totally dominate $c_1$, we
therefore have that $b_1 \in D$ or $d_1 \in D$.
If $v_2 \in D$, we let $D' = D \setminus \{v_1,v_2\}$. If $v_2 \notin
D$, then in order to totally dominate $v_1$ and $v_3$, we have that
$a_1 \in D$ and $c_1 \in D$, respectively. In this case, we let $D' =
D \setminus \{c_1,v_1\}$. In both cases, $D'$ is a TD-set of $G_k$
and $|D'| = |D| - 2$. Thus, $\gt(G_k) \le |D'| = \gt(F_k) - 2$.
Consequently, $\gt(F_k) = \gt(G_k) + 2$. Since $G_k \in \cG_\cub$ has
order~$4k$, we have by Theorem~B(c) that $\gt(G_k) = 2k$. Hence,
$\gt(F_k) = 2(k+1) = (n+1)/2$.

(b) If $G \ne C_3$ and $uv_1v_2v_3v$ is a special $2$-path in $G$,
then renaming $u$ and $v$, if necessary, we may assume $u = a_1$ and
$v = c_1$. Let $w \in \{v_1,v_2,v_3\}$. If $w = v_1$, let $S = C \cup
D \cup \{v_3\}$. If $w = v_2$, let $S = A \cup B \cup \{c_1\}$. If $w
= v_3$, let $S = A \cup B \cup \{v_1\}$. In all three cases, the set
$S$ is a TD-set of $G - w$ and $|S| = (n-1)/2$, and so $\gt(G - w)
\le |S| = (n-1)/2$. Every $\gt(G - w)$-set can be extended to a
TD-set of $G$ by adding to it a neighbor of $w$ in $G$, and so
$(n+1)/2 = \gt(G) \le \gt(G - w) + 1 \le (n+1)/2$. Consequently, we
must have equality throughout this inequality chain, implying that
$\gt(G - w) = (n-1)/2$.~

(c) We adopt the notation as in Part~(b) above. If $w = v_1$, let $S
= C \cup D \cup \{w\}$. If $w = v_2$ or $w = v_3$, let $S = A \cup B
\cup \{w\}$. In both cases the set $S$ is an ATD-set of $G$ with
respect to $w$ and $|S| = (n-1)/2$, and so $\gta(G;w) \le |S| =
(n-1)/2$. Every $\gta(G;w)$-set can be extended to a TD-set of $G$ by
adding to it a neighbor of $w$ in $G$, and so $(n+1)/2 = \gt(G) \le
\gta(G;w) + 1 \le (n+1)/2$. Consequently, we must have equality
throughout this inequality chain, implying that $\gta(G;w) =
(n-1)/2$.~\qed

\medskip
We note that if $G \in \cL$ and $G \ne C_6$, then by construction
there are only two $2$-paths in $G$ and both these $2$-paths are
special $2$-paths.

\begin{prop}
Let $G \in \cL$ have order~$n$. If $G \ne C_6$, let $uv_1v_2v_3v$ and
$u'u_1u_2u_3v'$ be the two special $2$-paths in $G$ and let $w \in
\{v_1,v_2,v_3\}$. Then the following holds. \\
\indent {\rm (a)} $\gt(G) = (n+2)/2$. \\
\indent {\rm (b)} $\gt(G - w) = n/2$. \\
\indent {\rm (c)} $\gta(G;w) = n/2$. \\
  \label{p:cL}
\end{prop}
\textbf{Proof.} Suppose that $G \in \cL$ has order~$n$. Then, $G =
L_k$ for some integer $k \ge 0$ and $n = 4k+6$. If $k = 0$, then $G =
C_6$. If $k \ge 1$, then following the notation introduced earlier,
the graph $L_k$ is obtained from the graph $G_k$ defined in
Section~\ref{S:GH} by subdividing the edge $a_1d_1$ three times and
subdividing the edge $b_kd_k$ three times. Equivalently, the graph
$L_k$ is obtained from the graph $F_k$ by subdividing the edge
$b_kc_k$ three times. Let $b_ku_1u_2u_3d_k$ denote the resulting path
in $L_k$. We note that the paths $a_1v_1v_2v_3c_1$ and
$b_ku_1u_2u_3d_k$ are the only $2$-paths and the only special
$2$-paths in~$G$.

(a) An analogous argument to show that $\gt(F_k) = \gt(G_k) + 2$ in
the proof of Proposition~\ref{p:cF} shows that $\gt(L_k) = \gt(F_k) +
2$. Hence, $\gt(L_k) = 2k + 4 = (n+2)/2$, as claimed.

(b) If $G \ne C_6$ and $uv_1v_2v_3v$ is a special $2$-path in $G$,
then renaming $u$ and $v$, if necessary, we may assume, by symmetry,
that $u = a_1$ and $v = c_1$. Let $w \in \{v_1,v_2,v_3\}$. If $w =
v_1$, let $S = C \cup D \cup \{u_1,u_2,v_3\}$. If $w = v_2$, let $S =
A \cup B \cup \{u_1,u_2,c_1\}$. If $w = v_3$, let $S = A \cup B \cup
\{u_1,u_2,v_1\}$. In all three cases, the set $S$ is a TD-set of $G -
w$ and $|S| = 2k + 3 = n/2$, and so $\gt(G - w) \le |S| = n/2$. Every
$\gt(G - w)$-set can be extended to a TD-set of $G$ by adding to it a
neighbor of $w$ in $G$, and so $(n+2)/2 = \gt(G) \le \gt(G - w) + 1
\le (n+2)/2$. Consequently, we must have equality throughout this
inequality chain, implying that $\gt(G - w) = n/2$.
%

(c) We adopt the notation as in Part~(b) above. If $w = v_1$, let $S
= C \cup D \cup \{u_1,u_2,w\}$. If $w = v_2$ or $w = v_3$, let $S = A
\cup B \cup \{u_1,u_2,w\}$. In both cases the set $S$ is an ATD-set
of $G$ with respect to $w$ and $|S| = n/2$, and so $\gta(G;w) \le |S|
= n/2$. Every $\gta(G;w)$-set can be extended to a TD-set of $G$ by
adding to it a neighbor of $w$ in $G$, and so $(n+2)/2 = \gt(G) \le
\gta(G;w) + 1 \le (n+2)/2$. Consequently, we must have equality
throughout this inequality chain, implying that $\gta(G;w) =
n/2$.~\qed

\begin{prop}
Let $G$ be a graph of order~$n$, size~$m$, total domination number
$\gt$, and maximum degree~$\Delta(G)$. Further let $\Delta = 3$ if
$\Delta(G) = 2$ and $\Delta = \Delta (G)$ if $\Delta(G) \ge 3$. If $G
\in \cG_{\done} \cup \cG_{\dtwo} \cup \cG_{\cub}$, then $m = \Delta
(n- \gt)$.
 \label{p:suff}
\end{prop}
\textbf{Proof.} Let $G \in \cG_{\done} \cup \cG_{\dtwo} \cup
\cG_{\cub}$ have order~$n$, size~$m$ and total domination number
$\gt$.
If $G \in \cG_{\done}$, then $G = C_k \circ P_2$ for some integer~$k
\ge 3$. Thus, $n = 3k$, $m = 3k$, and by Theorem~B(c), $\gt = 2k$,
implying that $m = 3(n - \gt)$.
If $G \in \cG_{\dtwo}$, then $G \in \cF \cup \cL$. If $G \in \cF$,
then $n = 4k+3$, $m = 6k + 3$, and by Proposition~\ref{p:cF}, $\gt =
2k+2$, implying that $m = 3(n - \gt)$. If $G \in \cL$, then $n =
4k+6$, $m = 6k + 6$, and by Proposition~\ref{p:cL}, $\gt = 2k+4$,
implying that $m = 3(n - \gt)$.
If $G \in \cG_{\cub}$, then $m = 3n/2$ and by Theorem~B(c), $\gt =
n/2$, implying that $m = 3(n - \gt)$.~\qed

\medskip
Following the notation introduced in Section~\ref{S:GH}
and~\ref{S:Gtwo}, we have the following useful property of graphs in
the family $\cG_{\dtwo} \cup \cG_{\cub}$.

\begin{prop}
Let $G \in \cG_{\dtwo} \cup \cG_{\cub}$ and let $v$ be an arbitrary
vertex in $G$. Then, $\gt(G - v) < \gt(G)$, unless one of the
following holds, in which case $\gta(G;v) < \gt(G)$. \\
\indent {\rm (a)} $G = G_k$ for some $k \ge 1$ and $v \in
\{a_1,b_k,c_1,d_k\}$. \\
\indent {\rm (b)} $G = F_k$ for some $k \ge 1$ and $v \in
\{b_k,d_k\}$.
 \label{p:cubic}
\end{prop}
\textbf{Proof.} Suppose first that $G \in \cG_{\cub}$. If $G =
GP_{16}$, then it is a simple exercise to check that $\gt(G - v) <
\gt(G)$. Suppose $G \in \cG$. Then, $G = G_k$ for some integer $k \ge
1$. Following the notation in Section~\ref{S:GH}, we may assume by
symmetry that $v = a_i$ for some $i$, $1 \le i \le k$. If $v = a_1$,
let $S_1 = (A \cup B) \setminus \{b_1\}$. Then, $S_1$ is an ATD-set
with respect to $v$, and so $\gta(G;v) \le |S_1| < |A| + |B| =
\gt(G)$. If $v = a_i$ where $i \ge 2$, let $S_i = (C \cup D)
\setminus \{d_i\}$. Then, $S_i$ is a TD-set in $G - v$, and so $\gt(G
- v) \le |S_i| < |C| + |D| = \gt(G)$.
Suppose $G \in \cH$. Then, $G = H_k$ for some integer $k \ge 2$.
Following the notation in Section~\ref{S:GH}, we may assume by
symmetry that $v = a_1$. Let $S = (A \cup B \cup \{c_1,d_k\})
\setminus \{a_1,b_1,b_k\}$. Then, $S$ is a TD-set in $G - v$, and so
$\gt(G - v) \le |S| < |A| + |B| = \gt(G)$.

Suppose next that $G \in \cG_{\dtwo}$. If $d(v) = 2$, then the result
follows from Propositions~\ref{p:cF} and~\ref{p:cL}. Hence we may
assume that $d(v) = 3$.

Suppose $G \in \cF$. Then, $G = F_k$ for some $k \ge 1$. Following
the notation introduced earlier, the graph $F_k$ is obtained from the
graph $G_k$ defined in Section~\ref{S:GH} by subdividing the edge
$a_1c_1$ three times. Let $a_1v_1v_2v_3c_1$ denote the resulting path
in $F_k$. We may assume by symmetry that $v = a_i$ or $v = b_i$ for
some $i$, $1 \le i \le k$.
%
%
If $v = a_i$ where $1 \le i \le k$, let $S_i' = (C \cup D \cup
\{v_2,v_3\}) \setminus \{d_i\}$. Then, $S_i'$ is a TD-set in $G - v$,
and so $\gt(G - v) \le |S_i'| < |C| + |D| + 2 = \gt(G)$.
If $v = b_k$, then let $D_k = (A \cup B \cup \{v_1,v_2\}) \setminus
\{a_k\}$. Then, $D_k$ is an ATD-set with respect to $v$, and so
$\gta(G;v) \le |D_k| < |A| + |B| + 2 = \gt(G)$.
If $v = b_i$ for some $i$, $1 \le i < k$, then let $D_i = (C \cup D
\cup \{v_2,v_3\}) \setminus \{c_i\}$. Then, $D_i$ is a TD-set in $G -
v$, and so $\gt(G - v) \le |D_i| < |C| + |D| + 2 = \gt(G)$.
Hence if $G \in \cF$, then the desired result follows.

Suppose $G \in \cL$. Then, $G = L_k$ for some $k \ge 1$. Following
the notation introduced earlier, the graph $L_k$ is obtained from the
graph $F_k$ by subdividing the edge $b_kd_k$ three times. Let
$b_ku_1u_2u_3d_k$ denote the resulting path in $L_k$. As observed
earlier, $d(v) = 3$. We may assume by symmetry that $v = a_i$ for
some $i$, $1 \le i \le k$. The set $(C \cup D  \cup
\{u_2,u_3,v_2,v_3\}) \setminus \{d_i\}$ is a TD-set in $G - v$, and
so $\gt(G - v) \le |C| + |D| + 3 < 2k + 4 = \gt(G)$. Hence if $G \in
\cF$, then $\gt(G - v) < \gt(G)$.~\qed

\section{Proof of Main Result}
\label{S:mainp}

Recall the statement of Theorem~\ref{thm1}.

\noindent \textbf{Theorem~\ref{thm1}}. \emph{Let $G$ be a connected
graph of order~$n$, size~$m$, total domination number $\gt$, and
maximum degree~$\Delta(G)$ with each component of $G$ of order at
least~$3$. Let $\Delta = 3$ if $\Delta(G) = 2$ and $\Delta = \Delta
(G)$ if $\Delta(G) \ge 3$. Then, $m \le \Delta (n- \gt)$, with
equality if and only if $G \in \cG_{\done} \cup \cG_{\dtwo} \cup
\cG_{\cub}$.}

\noindent \textbf{Proof of Theorem~\ref{thm1}.} The upper bound $m
\le \Delta (n- \gt)$ is a restatement of Theorem~A. If $G \in
\cG_{\done} \cup \cG_{\dtwo} \cup \cG_{\cub}$, then by
Proposition~\ref{p:suff}, $m = \Delta (n- \gt)$. Hence it suffices
for us to prove that if $m = \Delta (n- \gt)$, then $G \in
\cG_{\done} \cup \cG_{\dtwo} \cup \cG_{\cub}$. We proceed by
induction on the order~$n$ of $G$. If $G = P_3$, then $m = 2 < 3 =
\Delta (n- \gt)$. Hence if $n = 3$, then $G = C_3 \in \cG_{\dtwo}$.
This establishes the base case.

For the inductive hypothesis, let $n \ge 4$ and assume that if $G'$
is a connected graph of order~$n'$, size~$m'$, total domination
number $\gt'$ satisfying $m' = \Delta' (n'- \gt')$, where $\Delta' =
3$ if $\Delta(G') = 2$ and $\Delta' = \Delta(G')$ if $\Delta(G') \ge
3$, then $G' \in \cG_{\done} \cup \cG_{\dtwo} \cup \cG_{\cub}$. Let
$G$ be a connected graph of order~$n$, size~$m$, total domination
number $\gt$ satisfying $m = \Delta (n- \gt)$, where $\Delta = 3$ if
$\Delta(G) = 2$ and $\Delta = \Delta(G)$ if $\Delta(G) \ge 3$. We may
assume the following claim is satisfied by the graph $G$, for
otherwise the desired result holds.

\begin{unnumbered}{Claim~A.} The following hold in $G$. \\
\indent {\rm (a)} $\delta(G) \le 2$. \\
\indent {\rm (b)} $\Delta(G) \ge 3$.
\end{unnumbered}
\textbf{Proof.} (a) Suppose that $\delta(G) \ge 3$. By Theorem~B(b),
we have that $\gt \le n/2$. If $\gt < n/2$, then $m \le \Delta n / 2
= \Delta (n-n/2) < \Delta (n - \gt)$, a contradiction. Hence, $\gt =
n/2$, and so, by Theorem~B(c), $G \in \cG_{\cub}$, as desired.

(b) Suppose that $\Delta(G) = 2$. Then, $G$ is a path or a cycle. If
$G = P_n$, then $\gt \le 2n/3$, and so $m = n - 1 < 3(n - 2n/3) \le
\Delta (n - \gt)$, a contradiction. Hence, $G = C_n$. If $n \ne 6$,
then $\gt(G) < 2n/3$, and so $m = n = 3(n - 2n/3) < \Delta (n -
\gt)$, a contradiction. Therefore, $G = C_6 \in \cG_{\dtwo}$, as
desired.~\smallqed

\medskip
By Claim~A(b), $\Delta(G) \ge 3$, and so $\Delta = \Delta(G)$. Let
$\delta = \delta(G)$ and let $v$ be a vertex that has a neighbor, $w$
say, of degree~$\delta$. By Claim~A(a), $\delta \in \{1,2\}$. If
$V(G) = N[v]$, then $\Delta = n - 1 = d(v)$, $\gt = 2$ and $2m \le
\delta + (n-1) \Delta \le 2 + \Delta^2$. Thus since $\Delta \ge 3$,
we have $m \le (2 + \Delta^2)/2 < \Delta (\Delta - 1) = \Delta (n -
\gt)$, a contradiction. Hence, $V(G) \ne N[v]$.

Let $V_1$ be the set of isolated vertices in $G - N[v]$ and let $V_2$
be the set of vertices that belong to $P_2$-components of $G - N[v]$.
For $i = 1,2$, let $|V_i| = n_i$. If $V(G) \ne V_1 \cup V_2 \cup
N[v]$, let $F = G - N[v] - V_1 - V_2$ and let $H = G - V(F)$, i.e.,
$H = G[N[v] \cup V_1 \cup V_2]$. Then, $n = n(F) + n(H) = n(F) + d(v)
+ 1 + n_1 + n_2$ and $\gt \le \gt(F) + \gt(H)$. We proceed further
with the following claim.


\begin{unnumbered}{Claim~B.}
$n_1 + n_2 = 0$.
\end{unnumbered}
\textbf{Proof.} For the sake of contradiction, suppose that $n_1 +
n_2 \ge 1$. Then the following two claims established in~\cite{He05}
hold.

\begin{unnumbered}{Claim~B.1}
$\gt(H) \le n_1 + n_2 + 1$. Further if a vertex in $N(v)$ has two or
more neighbors in $V_1 \cup V_2$, then $\gt(H) \le n_1 + n_2$.
\end{unnumbered}

\begin{unnumbered}{Claim~B.2}
$\gt(H) \le d(v) + 1 + n_2/2$. Further if $n_1 = 0$, then $\gt(H) \le
d(v) + n_2/2$.
\end{unnumbered}

Let $m_1$ denote the number of edges of $G$ incident with vertices in
$N(v)$. Then, $m_1 \le \Delta (d(v) - 1) + \delta$. Since each
component of $F$ has order at least~$3$, applying Theorem~A to each
component of $F$ we deduce that $m(F) \le \Delta ( n(F) - \gt(F) )$.
Hence,

\[
\begin{array}{lclr} \2
\Delta ( n -\gt) & = & m & \\      \2
& \le & m_1 + m(G[V_2]) + m(F) & (1) \\     \2
& \le & m_1 + n_2/2 + \Delta ( n(F) - \gt(F) ) & (2) \\     \2
& \le & m_1 + n_2/2 + \Delta ( n - d(v) - 1 - n_1 - n_2 -
\gt + \gt(H) ) & (3) \\            \2
& = & \Delta ( n -\gt) + m_1 + n_2/2 - \Delta ( d(v) + 1 +
n_1 + n_2 ) + \Delta  \cdot \gt(H) &   \\            \2
& \le & \Delta ( n -\gt) - 2 \Delta + \delta + n_2/2 - \Delta
(n_1 + n_2 ) + \Delta \cdot \gt(H).   & (4)
\end{array}
\]
Let
\[
\xi(G) = - 2 \Delta + \delta + n_2/2 - \Delta (n_1 + n_2 ) +
\Delta \cdot \gt(H),
\]

\noindent and so by the above inequality chain, we have that $\xi(G)
\ge 0$.

\begin{unnumbered}{Claim~B.3}
$n_2 \le 2(\Delta -\delta)$. 
\end{unnumbered}
\textbf{Proof.} Assume, to the contrary, that $n_2 > 2(\Delta
-\delta)$. Then since $n_2$ is even, we have $n_2 \ge 2(\Delta
-\delta + 1)$.
Suppose first that $n_1 \ge 1$. By Claim~B.2, we have $\gt(H) \le
d(v) + n_2/2 + 1$. Hence if $n_2
> 2( \Delta \cdot d(v) + \delta - \Delta - \Delta n_1 )/ (\Delta -
1)$, then $\xi(G) < 0$, a contradiction. Therefore, $n_2 \le 2(
\Delta \cdot d(v) + \delta - \Delta - \Delta n_1 )/ (\Delta - 1) =
2d(v) + 2( d(v) + \delta - \Delta - \Delta n_1 )/ (\Delta - 1)$.
However since $d(v) \le \Delta$ and since, by Claim~A(a), $\delta \le
2 < \Delta$, we have that $d(v) + \delta < 2\Delta$. Further since
$n_1 \ge 1$, we have that $ - \Delta n_1 \le - \Delta$, implying that
$n_2 < 2 d(v)$.
Suppose next that $n_1 = 0$. By Claim~B.2, we have $\gt(H) \le d(v) +
n_2/2$. If $n_2 > 2( \Delta \cdot d(v) +  \delta - 2 \Delta )/
(\Delta - 1)$, then $\xi(G) < 0$, a contradiction. Therefore, $n_2
\le 2( \Delta \cdot d(v) +  \delta - 2 \Delta  )/ (\Delta - 1) =
2d(v) + 2( d(v) + \delta - 2 \Delta  )/ (\Delta - 1) < 2 d(v)$.
Thus irrespective of whether $n_1 \ge 1$ or $n_1 = 0$, we have
\[
\hspace*{5cm} 2( \Delta - \delta  + 1) \le n_2 < 2 d(v), \hspace*{5cm} (5)
\]
which implies that $d(v) \ge \Delta - \delta + 2$. Since $\Delta \ge
d(v)$, this in turn implies that $\delta \ge 2$. Consequently by
Claim~A(a) we have $\delta = 2$. Thus there are at least $n_2 \ge
2(\Delta - 1)$ edges between $V_2$ and $N(v)$.

Since $d(v) \le \Delta$, the inequality chain~(5)
implies that $n_2 < 2 \Delta$. Hence if $\gt(H) \le n_1 + n_2$, then
we have $\xi(G) \le - 2 \Delta + \delta + n_2/2 < - \Delta + \delta <
0$, a contradiction. Therefore, $\gt(H) \ge n_1 + n_2 + 1$.
Consequently, by Claim~B.1, $\gt(H) = n_1 + n_2 + 1$ and every vertex
in $N(v)$ has at most one neighbor in $V_1 \cup V_2$. Hence, $\Delta
\ge d(v) \ge 2( \Delta - 1)$, and so $\Delta \le 2$, a
contradiction.~\smallqed

\medskip
By Claim~B.3, we have $n_2 \le 2(\Delta -\delta)$. If $\gt(H) \le n_1
+ n_2$, then $\xi(G) < 0$, a contradiction. Hence, $\gt(H) \ge n_1 +
n_2 + 1$. Consequently by Claim~B.1, $\gt(H) = n_1 + n_2 + 1$, and so
$\xi(G) = - \Delta + \delta + n_2/2$. By Claim~B.1, we note that
every vertex in $N(v)$ has at most one neighbor in $V_1 \cup V_2$. If
$n_2 < 2(\Delta -\delta)$, then $\xi(G) < 0$, a contradiction. Hence,
$n_2 = 2(\Delta -\delta)$ and $\xi(G) = 0$. But this implies that we
must have equality throughout the inequality chain following
Claim~B.2 (and preceding Claim~B.3). In particular, equality in~(1),
(2), (3) and (4) implies the following claim.

\begin{unnumbered}{Claim~B.4} The following hold in $G$. \\
\indent {\rm (a)} $N(v)$ is an independent set. \\
\indent {\rm (b)} $m(F) = \Delta ( n(F) - \gt(F) )$. \\
\indent {\rm (c)} $m(F') = \Delta ( n(F') - \gt(F') )$ for each
component $F'$ of $F$. \\
\indent {\rm (d)} $\gt = \gt(F) + \gt(H)$. \\
\indent {\rm (e)} $m_1 = \Delta (d(v) - 1) + \delta$. \\
\indent {\rm (f)} Every neighbor of $v$ different from $w$ has
degree~$\Delta$.
\end{unnumbered}

As observed earlier, every vertex in $N(v)$ has at most one neighbor
in $V_1 \cup V_2$. Hence since $N(v)$ is an independent set, every
neighbor of $v$ of degree~$\Delta$ has at least~$\Delta - 2 \ge 1$
neighbors in $V(F)$. Applying the inductive hypothesis to every
component $F'$ of $F$, we have $F' \in \cG_{\done} \cup \cG_{\dtwo}
\cup \cG_{\cub}$. In particular, we note that $\Delta(F) = 3$. Thus
by Theorem~A, we have $m(F) \le 3(n(F) - \gt(F))$, and so by
Claim~B.4(b) we have $\Delta \le 3$. Consequently, $\Delta = 3$.

Suppose $F' \in \cG_{\cub}$ for some component $F'$ of $F$. Since $G$
is connected, there is a vertex $v'$ in $F'$ adjacent to a neighbor
of $v$. But then $\Delta \ge d_G(v') \ge d_F(v') + 1 = 4$, a
contradiction. Hence, $F' \in \cG_{\done} \cup \cG_{\dtwo}$ for every
component $F'$ of $F$.

Let $N_v = N(v) \setminus \{w\}$. If no vertex in $V_1 \cup V_2$ has
a neighbor in $N_v$, then since $G$ is connected, the vertex $w$ is
adjacent to every vertex in $V_1 \cup V_2$. However, $d_G(w) = \delta
\le 2$, implying that $\delta = 2$ and $n_1 + n_2 = 1$. Let $N(w) =
\{v,w'\}$. On the one hand, if $w' \in V_1$, then since $w'$ has
degree at least~$2$ in $G$, the vertex $w'$ has a neighbor in $N_v$.
On the other hand, if $w' \in V_2$, then the neighbor of $w'$ in
$G[V_2]$ has a neighbor in $N_v$. Both cases produce a contradiction.
Hence there is a vertex $u' \in V_1 \cup V_2$ that is adjacent to a
vertex $u \in N_v$. Since every vertex in $N(v)$ has at most one
neighbor in $V_1 \cup V_2$, and since $d_G(u) = \Delta = 3$, the
vertex $u$ is therefore adjacent in $G$ to a vertex $z \in V(F)$.
Since $\Delta = 3$, the vertex $z$ has either degree~$1$ in $F$ or
degree~$2$ in $F$.

We now construct a $\gt(H)$-set $S_H$ as follows. Initially, let $S_H
= \{u,v\}$. On the one hand, suppose $u' \in V_1$. For every vertex
in $V_1 \setminus \{u'\}$, choose a neighbor in $N(v)$ and add it to
$S_H$. Further for every $K_2$-component in $G[V_2]$, choose a vertex
that has a neighbor in $N(v)$ and add both the chosen vertex in $V_2$
and one of its neighbors in $N(v)$ to the set $S_H$. On the other
hand, suppose $u' \in V_2$. Then add $u'$ to $S_H$ and for every
$K_2$-component in $G[V_2]$ that does not contain $u'$, choose a
vertex that has a neighbor in $N(v)$ and add both the chosen vertex
in $V_2$ and one of its neighbors in $N(v)$ to the set $S_H$. Further
for every vertex in $V_1$, choose a neighbor in $N(v)$ and add it to
$S_H$. In both cases, the resulting set $S_H$ is a TD-set of $H$ and
$|S_H| = n_1 + n_2 + 1$. Thus, $S_H$ is a $\gt(H)$-set that contains
the vertex $u$.

If $F' \in \cG_{\done}$ and $d_F(z) = 1$, let $S_F'$ be a $\gt(F' -
z)$-set. If $F' \in \cG_{\done}$ and $d_F(z) = 2$, let $S_F'$ be a
$\gta(F';z)$-set. If $F' \in \cG_{\dtwo}$, then $d_F(z) = 2$ and
necessarily $z$ is an internal vertex of a special $2$-path in $F'$,
and we let $S_F'$ be a $\gt(F' - z)$-set. By
Observation~\ref{ob:Gone}, and by Propositions~\ref{p:cF}
and~\ref{p:cL} we have that $|S_F'| = \gt(F') - 1$. If $F = F'$, we
let $S_F = S_F'$. If $F$ contains at least two components, then we
add to $S_F'$ a $\gt(F^*)$-set from every component $F^*$ of $F$
different from $F'$ and we let $S_F$ denote the resulting set. In
both cases, $|S_F| = \gt(F) - 1$. By construction the set $S_F \cup
S_H$ is a TD-set of $G$, and so $\gt \le |S_F \cup S_H| = \gt(H) +
\gt(F) - 1$, contradicting Claim~B.4(d). Therefore, $n_1 + n_2 = 0$.
This completes the proof of Claim~B.~\smallqed

\medskip
By Claim~B, $n_1 + n_2 = 0$, and so $H = G[N[v]]$ and $n = n(H) +
n(F)$. Further, $\gt \le \gt(H) + \gt(F) = \gt(F) + 2$. Recall that
$v$ is a vertex in $G$ with a neighbor, $w$, of degree~$\delta$.
Recall also that $\delta \in \{1,2\}$ and $\Delta(G) = \Delta \ge 3$.

\newpage
\begin{unnumbered}{Claim~C.}
If $\delta = 1$, then every support vertex has degree~$2$ in $G$.
\end{unnumbered}
\textbf{Proof.} Let $\delta = 1$. For the sake of contradiction,
assume that there is a support vertex in $G$ whose degree is at
least~$3$. Renaming vertices, if necessary, we may assume that
$d_G(v) \ge 3$. Let $x$ be an arbitrary vertex in $N(v) \setminus
\{w\}$ and let $G_x = G - (N(v) \setminus \{w,x\})$. In particular,
we note that $V(F) \subset V(G_x)$ and that $F$ is an induced
subgraph in $G_x$. Every TD-set of $G_x$ must contain the support
vertex $v$, and is therefore also a TD-set of $G$. Hence, $\gt =
\gt(G) \le \gt(G_x)$. Since each component of $F$, and therefore of
$G_x$, has order at least~$3$, applying Theorem~A to $G_x$ we obtain
$m(G_x) \le \Delta ( n(G_x) - \gt(G_x) )$. Therefore,

\[
\begin{array}{lcl} \1
\Delta ( n - \gt) & = & m \\          \1
 & \le & m(G_x) + \Delta \, (d(v) - 2)  \\    \1
& \le & \Delta ( n(G_x) - \gt(G_x) ) + \Delta \, (d(v) - 2) \\      \1
& \le & \Delta ( n - d(v) + 2 - \gt ) + \Delta \, (d(v) - 2) \\     \1
& = & \Delta ( n -\gt).
\end{array}
\]

Hence we must have equality throughout the above inequality chain.
Since $x$ is an arbitrary vertex in $N(v) \setminus \{w\}$, this
implies that $N(v)$ is an independent set and every neighbor of $v$
different from $w$ has degree~$\Delta$. Further, $m(G_x) = \Delta (
n(G_x) - \gt(G_x) )$, and so $m(F_x) = \Delta ( n(F_x) - \gt(F_x) )$
for each component $F_x$ of $G_x$. Applying the inductive hypothesis
to every component $F_x$ of $G_x$, we have $F_x \in \cG_{\done} \cup
\cG_{\dtwo} \cup \cG_{\cub}$. In particular, we note that
$\Delta(G_x) = 3$. Thus by Theorem~A, we have $m(G_x) \le 3(n(G_x) -
\gt(G_x))$, implying that $\Delta \le 3$. Consequently, $\Delta = 3$.
In particular, $d(v) = 3$. Let $N(v) = \{w,x,y\}$.

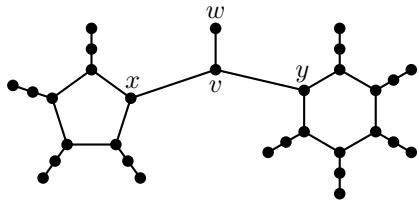
\begin{figure}[htb]
\tikzstyle{every node}=[circle, draw, fill=black!100, inner sep=0pt,minimum width=.125cm]
\begin{center}
\begin{tikzpicture}[thick,scale=.55]
  \draw(0,0) { 
    +(4.90,4.00) -- +(4.90,3.00)
    +(4.90,3.00) -- +(2.85,2.31)
    +(4.90,3.00) -- +(7.04,2.50)
    +(7.04,2.50) -- +(7.90,3.00)
    +(7.90,3.00) -- +(8.77,2.50)
    +(8.77,2.50) -- +(8.77,1.50)
    +(8.77,1.50) -- +(7.90,1.00)
    +(7.90,1.00) -- +(7.04,1.50)
    +(7.04,1.50) -- +(7.04,2.50)
    +(7.90,3.00) -- +(7.90,3.50)
    +(7.90,3.50) -- +(7.90,4.00)
    +(8.77,2.50) -- +(9.20,2.75)
    +(9.20,2.75) -- +(9.63,3.00)
    +(8.77,1.50) -- +(9.20,1.25)
    +(9.20,1.25) -- +(9.63,1.00)
    +(7.90,1.00) -- +(7.90,0.50)
    +(7.90,0.50) -- +(7.90,0.00)
    +(7.04,1.50) -- +(6.60,1.25)
    +(6.60,1.25) -- +(6.17,1.00)
    +(2.85,2.31) -- +(1.90,3.00)
    +(1.90,3.00) -- +(0.95,2.31)
    +(0.95,2.31) -- +(1.31,1.19)
    +(1.31,1.19) -- +(2.49,1.19)
    +(2.49,1.19) -- +(2.85,2.31)
    +(2.49,1.19) -- +(2.78,0.79)
    +(2.78,0.79) -- +(3.08,0.38)
    +(1.31,1.19) -- +(1.02,0.79)
    +(1.02,0.79) -- +(0.73,0.38)
    +(0.95,2.31) -- +(0.48,2.46)
    +(0.48,2.46) -- +(0.00,2.62)
    +(1.90,4.00) -- +(1.90,3.50)
    +(1.90,3.50) -- +(1.90,3.00)
    +(4.90,3.00) node{}
    +(4.90,4.00) node{}
    +(0.95,2.31) node{}
    +(1.90,3.00) node{}
    +(2.85,2.31) node{}
    +(2.49,1.19) node{}
    +(1.31,1.19) node{}
    +(0.48,2.46) node{}
    +(0.00,2.62) node{}
    +(1.02,0.79) node{}
    +(0.73,0.38) node{}
    +(2.78,0.79) node{}
    +(3.08,0.38) node{}
    +(1.90,3.50) node{}
    +(1.90,4.00) node{}
    +(7.04,2.50) node{}
    +(7.90,3.00) node{}
    +(8.77,2.50) node{}
    +(8.77,1.50) node{}
    +(7.90,1.00) node{}
    +(7.04,1.50) node{}
    +(6.60,1.25) node{}
    +(6.17,1.00) node{}
    +(9.20,1.25) node{}
    +(9.63,1.00) node{}
    +(7.90,0.50) node{}
    +(7.90,0.00) node{}
    +(9.20,2.75) node{}
    +(9.63,3.00) node{}
    +(7.90,3.50) node{}
    +(7.90,4.00) node{}
    +(4.9,2.6) node[rectangle, draw=white!0, fill=white!100]{$v$}
    +(4.9,4.4) node[rectangle, draw=white!0, fill=white!100]{$w$}
    +(2.9,2.7) node[rectangle, draw=white!0, fill=white!100]{$x$}
    +(7,2.9) node[rectangle, draw=white!0, fill=white!100]{$y$} 
  };
\end{tikzpicture}
\end{center}
\vskip -0.6 cm \caption{The graph $G$.} \label{FxFy}
\end{figure}

Let $F_x$ be the component of $G_x$ that contains $x$ and let $F_y$
be the component of $G_y$ that contains~$y$. Since $F_x$ contains a
vertex of degree~$1$, applying the inductive hypothesis to $F_x$ we
have that $F_x \in \cG_{\done}$. Analogously, $F_y \in \cG_{\done}$.
Further, $F_x$ and $F_y$ have only the vertices $v$ and $w$ in
common. Since $G$ is connected, we note that $V(G) = V(F_x) \cup
V(F_y)$. The component of $G - vx$ containing $v$ is $F_y$ and the
component of $G - vy$ containing $v$ is $F_x$. For example, if the
cycle in $F_x$ has length~$5$ and the cycle in $F_y$ has length~$6$,
then the graph $G$ is illustrated in Figure~\ref{FxFy}. If the cycle
in $F_x$ has length~$k_1$ and the cycle in $F_y$ has length~$k_2$,
then $n = 3k_1 + 3k_2 - 2$, $m = 3k_1 + 3k_2 - 1$, and $\gt = 2(k_1 -
1) + 2(k_2 - 1) + 2 = 2(k_1 + k_2 - 1)$. But then $m < \Delta ( n -
\gt )$, a contradiction.~\smallqed

By Claim~C, every support vertex has degree~$2$ in $G$. In
particular, $d(v) = 2$. Let $N(v) = \{u,w\}$. If $n = \Delta + 2$,
then $\gt = 2$ and $m \le {\Delta \choose 2} + 2 < \Delta^2 =
\Delta(n - \gt)$, a contradiction. Hence, $n \ge \Delta + 3$.

\begin{unnumbered}{Claim~D.}
If $\delta = 1$, then the vertex at distance~$2$ from a vertex of
degree~$1$ has degree~$\Delta$ in $G$.
\end{unnumbered}
\textbf{Proof.} For the sake of contradiction, assume that there is a
vertex at distance~$2$ from a vertex of degree~$1$ with degree less
than~$\Delta$ in $G$. Renaming vertices, if necessary, we may assume
that $d(u) < \Delta$. Then,

\[
\begin{array}{lcl} \1
\Delta ( n - \gt) & = & m \\          \1
 & = & d(u) + d(w) + m(F)  \\    \1
& \le & \Delta + \Delta ( n(F) - \gt(F) )  \\      \1
& \le & \Delta + \Delta ( n - 3 - \gt + 2 ) \\     \1
& = & \Delta (n -\gt).
\end{array}
\]

Hence we must have equality throughout the above inequality chain. In
particular, $\gt = \gt(F) + 2$ and $m(F) = \Delta ( n(F) - \gt(F) )$,
and so $m(F') = \Delta ( n(F') - \gt(F') )$ for each component $F'$
of $F$. Applying the inductive hypothesis to every component $F'$ of
$F$, we have $F' \in \cG_{\done} \cup \cG_{\dtwo} \cup \cG_{\cub}$,
implying that $\Delta = 3$ and $F' \in \cG_{\done} \cup \cG_{\dtwo}$.
Let $z \in N(u) \setminus \{v\}$ and let $F_z$ be the component of
$F$ containing $z$ (possibly, $F_z = F$). Since $\Delta = 3$, the
vertex $z$ has either degree~$1$ in $F$ or degree~$2$ in $F$. Let
$S_H = \{u,v\}$. Proceeding as in the last paragraph of the proof of
Claim~B, the set $S_H$ can be extended to a TD-set of $G$ of
cardinality at most~$\gt(H) + \gt(F) - 1 = \gt(F) + 1$. Thus, $\gt
\le \gt(F) + 1$, contradicting our earlier observation that $\gt =
\gt(F) + 2$.~\smallqed

\begin{unnumbered}{Claim~F.}
If $\delta = 1$, then $G \in \cG_{\done}$.
\end{unnumbered}
\textbf{Proof.} Suppose $\delta = 1$. By Claim~C, $d(v) = 2$ and by
Claim~D, $d(u) = \Delta$. We proceed further with the following
claim.

\begin{unnumbered}{Claim~F.1.}
If $N(u) \setminus \{v\}$ does not induce a clique, then $G \in
\cG_{\done}$.
\end{unnumbered}
\textbf{Proof.} Suppose that $N(u) \setminus \{v\}$ does not induce a
clique. Let $x$ and $y$ be two vertices in $N(u) \setminus \{v\}$
that are not adjacent and let $L$ be obtained from $F$ by adding to
it the edge $xy$. Each vertex in $N(u) \setminus \{v\}$ has degree at
most~$\Delta - 1$ in $F$, and therefore degree at most~$\Delta$ in
$L$, implying that $\Delta(L) \le \Delta$.  Every TD-set of $L$ can
be extended to a TD-set of $G$ by adding to it the vertices $u$ and
$v$, and so $\gt \le \gt(L) + 2$. Since each component of $L$ has
order at least~$3$, applying Theorem~A to $L$ we have that $m(L) \le
\Delta( n(L) - \gt(L) )$. Thus,
\[
\begin{array}{lcl} \1
\Delta ( n - \gt) & = & m \\          \1
 & = & 1 + \Delta + (m(L) - 1) \\    \1
& \le & \Delta + \Delta ( n(L) - \gt(L) )  \\      \1
& \le & \Delta + \Delta ( n - 3 - \gt + 2 ) \\     \1
& = & \Delta (n -\gt).
\end{array}
\]

Hence we must have equality throughout the above inequality chain. In
particular, $\gt = \gt(L) + 2$ and $m(L) = \Delta ( n(L) - \gt(L) )$,
and so $m(L') = \Delta ( n(L') - \gt(L') )$ for each component $L'$
of $L$. Applying the inductive hypothesis to every component $L'$ of
$L$, we have $L' \in \cG_{\done} \cup \cG_{\dtwo} \cup \cG_{\cub}$,
implying that $\Delta = 3$. In particular, $N(u) \setminus \{v\} =
\{x,y\}$, and so $L$ is connected since $G$ is connected.

Since $\delta = 1$, both $x$ and $y$ have degree at least~$1$ in $F$
and therefore at least~$2$ in $L$. Suppose $x$ or $y$, say $x$, has
degree~$2$ in $L$. Then, $L \in \cG_{\dtwo}$ and $x$ is an internal
vertex of a special $2$-path in $L$. Let $S_H = \{u,v\}$. Proceeding
as in the last paragraph of the proof of Claim~B, the set $S_H$ can
be extended to a TD-set of $G$ of cardinality at most~$\gt(H) +
\gt(L) - 1 = \gt(L) + 1$. Thus, $\gt \le \gt(L) + 1$, contradicting
our earlier observation that $\gt = \gt(L) + 2$. Hence both $x$ and
$y$ have degree~$3$ in $L$.

If $L \in \cG_{\dtwo} \cup \cG_{\cub}$, then by
Proposition~\ref{p:cubic}, either $\gt(L - x) < \gt(L)$ or $\gta(L;x)
< \gt(L)$. However every $\gt(L - x)$-set and every $\gta(L;x)$-set
can be extended to a TD-set of $G$ by adding to it the set $\{u,v\}$,
implying that $\gt < \gt(L) + 2$, a contradiction. Hence, $L \in
\cG_{\done}$. But then $G \in \cG_{\done}$.~\smallqed

\medskip
By Claim~F.1, we may assume that $N(u) \setminus \{v\}$ induces a
clique. Thus each vertex in $N(u) \setminus \{v\}$ has degree at
least~$\Delta - 1$. More generally, we may assume with our
assumptions to date that the neighbor of every vertex $v_1$ of
degree~$1$ is a vertex $v_2$ of degree~$2$ whose other neighbor $v_3$
is a vertex of degree~$\Delta$, and that $N(v_3) \setminus \{v_2\}$
induces a clique. Let $x \in N(u) \setminus \{v\}$ and let $F_x = G -
\{u,v,w,x\}$.

\begin{unnumbered}{Claim~F.2.}
The graph $F_x$ has an isolated vertex or a $P_2$-component.
\end{unnumbered}
\textbf{Proof.} Suppose that $F_x$ has no isolated vertex and no
$P_2$-component. Applying Theorem~A to $F_x$, we obtain

\[
\begin{array}{lcl}
m & \le & 2 \Delta + m(F_x)\\
& \le & 2 \Delta + \Delta(n(F_x) - \gt(F_x))\\
 &\le & 2 \Delta + \Delta(n - 4 + 2 -\gt)\\
 & = &\Delta(n-\gt).
\end{array}
\]

Hence we must have equality throughout the above inequality chain. In
particular, $d_G(x) = \Delta$, $\gt = \gt(F_x) + 2$ and $m(F_x) =
\Delta ( n(F_x) - \gt(F_x) )$, and so $m(F') = \Delta ( n(F') -
\gt(F') )$ for each component $F'$ of $F_x$. Applying the inductive
hypothesis to every component $F'$ of $F_x$, we have $F' \in
\cG_{\done} \cup \cG_{\dtwo} \cup \cG_{\cub}$, implying that $\Delta
= 3$ and that $F' \in \cG_{\done} \cup \cG_{\dtwo}$. Let $N(u) =
\{u,x,y\}$ and consider the component $F'$ of $F_x$ that
contains~$y$. Since $\Delta = 3$ and $y$ is adjacent to both $u$ and
$x$ in $G$, we have that $F' \in \cG_{\done}$ and $y$ has degree~$1$
in $F'$. A $\gt(F' - y)$-set can be extended to a $\gt(G)$-set of
cardinality less than~$\gt$ by adding to it the set $\{v,w\}$ and, if
$F_x$ is not connected, adding a minimum TD-set from each component
of $F_x$ not containing~$y$, a contradiction.~\smallqed

\medskip
By Claim~F.2, for every $x \in N(u) \setminus \{v\}$, the graph $F_x$
has an isolated vertex or a $P_2$-component. Recall that $n \ge
\Delta + 3$. Since $F$ has no isolated vertex and no $P_2$-component,
the vertex $x$ must be adjacent to each isolated vertex of $F_x$ and
to at least one vertex from each $P_2$-component of $F_x$. Thus it
follows from our earlier assumptions that $d(x) = \Delta$ and that
$F_x$ has no isolated vertex and exactly one $P_2$-component that is
joined to $x$ by exactly one edge. Hence, $G = K_{\Delta} \circ P_2$,
and so $n = 3\Delta$, $\gt = 2\Delta$, and $m = (\Delta^2 + 3
\Delta)/2$. If $\Delta > 3$, then $m < \Delta (n - \gt)$, a
contradiction. Hence, $\Delta = 3$ and $G = K_{\Delta} \circ P_2 \in
\cG_{\done}$. This completes the proof of Claim~F.~\qed

\medskip
By Claim~F, we may assume that $\delta =2$. Let $S_2$ denote the set
of vertices of $G$ that have degree~$2$. If every component of
$G[S_2]$ has order at most~$2$, then by Theorem~B(d), $\gt \le n/2$.
Since $\delta = 2$ and $\Delta = \Delta(G) \ge 3$, we have that $m <
\frac{1}{2} \Delta \, n= \Delta (n-\frac{1}{2}n) \le \Delta (n-\gt)$,
a contradiction. Hence, $G[S_2]$ has a path of length at least~$2$.
Renaming vertices if necessary, we may assume that $vwx$ is such a
path in $G[S_2]$. Let $N(v) = \{u,w\}$ and let $N(x) = \{w,y\}$.

By our assumptions to date, if $z$ is a vertex that is adjacent to a
vertex of degree~$2$ in $G$, then every component of $G - N[z]$ has
order at least~$3$. In particular, since $u$ is adjacent to the
vertex $v$ of degree~$2$ in $G$, every component of $G - N[u]$ has
order at least~$3$. If $u = y$, then $w$ is an isolated vertex in $G
- N[u]$, a contradiction. Hence, $u \ne y$. If $u$ and $y$ are
adjacent, then $wx$ is a $P_2$-component in $G - N[u]$, a
contradiction. Hence, $u$ and $y$ are not adjacent.

Let $L$ be obtained from $G - \{v,w,x\}$ by adding the edge $uy$.
Since $G$ is connected, so too is $L$. Each of $u$ and $y$ has degree
at most~$\Delta$ in $L$ and degree at least~$2$ in $L$, implying that
$2 \le \delta(L)$ and $\Delta(L) \le \Delta$. Applying Theorem~A to
$L$, we have $m(L) \le \Delta (n(L) - \gt(L))$. Let $S_L$ be a
$\gt(L)$-set. If $\{u,y\} \subseteq L$, let $S = S_L \cup \{v,x\}$.
If $u \in S_L$ and $y \notin S_L$, let $S = S_L \cup \{w,x\}$. If $u
\notin S_L$, let $S = S_L \cup \{v,w\}$. In all three cases, $S$ is a
TD-set of $G$, and so $\gt \le |S| \le |S_L| + 2 = \gt(L) + 2$.
Counting edges in $G$, we therefore have that
\[
\begin{array}{lcl}
m & = & 4 + (m(L) - 1)\\
& \le & 3 + \Delta (n(L) - \gt(L))\\
 &\le & 3 + \Delta(n - 3 - \gt + 2)\\
 & = & \Delta(n-\gt) - \Delta + 3 \\
  & \le & \Delta(n-\gt).
\end{array}
\]

Hence we must have equality throughout the above inequality chain. In
particular, $\Delta = 3$, $\gt = \gt(L) + 2$ and $m(L) = \Delta (n(L)
- \gt(L))$. Applying the inductive hypothesis to $L$, we have $L \in
\cG_{\dtwo} \cup \cG_{\cub}$.

\begin{unnumbered}{Claim~G.}
If $L \in \cG_{\cub}$, then $G \in \cG_{\dtwo}$.
\end{unnumbered}
\textbf{Proof.} Suppose $L \in \cG_{\cub}$. If $\gt(L - u) < \gt(L)$,
then since every $\gt(L - u)$-set can be extended to a TD-set in $G$
by adding to it the vertices $v$ and $w$, we have that $\gt < \gt(L)
+ 2$, a contradiction. Hence, $\gt(L - u) = \gt(L)$. Analogously,
$\gt(L - y) = \gt(L)$. Therefore by Proposition~\ref{p:cubic}, we
have that $L = G_k \in \cG$ for some $k \ge 1$ and that $\{u,y\}
\subset \{a_1,b_k,c_1,d_k\}$ (following the notation in
Section~\ref{S:GH}). Renaming vertices, if necessary, we may assume
that $\{u,y\} = \{a_1,c_1\}$. But then $G = F_k \in
\cG_{\dtwo}$.~\smallqed

\begin{unnumbered}{Claim~H.}
If $L \in \cG_{\dtwo}$, then $G \in \cG_{\dtwo}$.
\end{unnumbered}
\textbf{Proof.} Suppose $L \in \cG_{\dtwo}$. If $\gt(L - u) <
\gt(L)$, then $\gt < \gt(L) + 2$, a contradiction. Hence, $\gt(L - u)
= \gt(L)$. Analogously, $\gt(L - y) = \gt(L)$. Therefore by
Proposition~\ref{p:cubic}, we have that $L = F_k \in \cF$ for some $k
\ge 1$ and that $\{u,y\} = \{b_k,d_k\}$. But then $G = L_k \in \cL
\subset \cG_{\dtwo}$.~\smallqed

\medskip
Hence if $\delta = 2$, then by Claim~G and Claim~H, we have $G \in
\cG_{\dtwo}$. This completes the proof of Theorem~\ref{thm1}.~\qed

\end{document}